\documentclass[final,1p,times]{elsarticle}
\usepackage{amsmath}
\usepackage{amssymb}
\usepackage{url}
\usepackage{amsfonts}
\usepackage{graphicx}
\usepackage{rotating}
\usepackage{floatflt,epsfig}
\usepackage{lineno,hyperref}
\usepackage{enumerate}
\usepackage{colortbl}
\usepackage{array,tabularx,tabulary,booktabs}
\usepackage{longtable}
\usepackage{multirow}
\usepackage{wrapfig}
\usepackage{subcaption}
\newcolumntype{^}{>{\currentrowstyle}}

\makeatletter
\def\ps@pprintTitle{%
   \let\@oddhead\@empty
   \let\@evenhead\@empty
   \let\@oddfoot\@empty
   \let\@evenfoot\@oddfoot
}
\makeatother

\newtheorem{lemma}{Lemma}

\newtheorem{theorem}{Theorem}
\newtheorem{corollary}{Corollary}
\newtheorem{proposition}{Proposition}

\newtheorem{remark}{Remark}
\newtheorem{construction}{Construction}

\newcommand{\proof}{\medskip\noindent{\bf Proof.~}}
\begin{document}
\renewcommand{\abstractname}{Abstract}
\renewcommand{\refname}{References}
\renewcommand{\arraystretch}{0.9}
\thispagestyle{empty}
\sloppy

\begin{frontmatter}
\title{Divisible design graphs with parameters $(4n,n+2,n-2,2,4,n)$\\
and $(4n,3n-2,3n-6,2n-2,4,n)$}

\author[01]{Leonid~Shalaginov}
\ead{44sh@mail.ru}

\address[01]{Chelyabinsk State University, Brat'ev Kashirinyh st. 129, Chelyabinsk, 454021, Russia}

\begin{abstract}
A $k$-regular graph is called a divisible design graph (DDG for short) if its vertex set can be partitioned into $m$ classes of size $n$, such that two distinct vertices from the same class have exactly $\lambda_1$ common neighbors, and two vertices from different classes have exactly $\lambda_2$ common neighbors. $4\times n$-lattice graph is the line graph of $K_{4,n}$. This graph is a DDG with parameters $(4n,n+2,n-2,2,4,n)$. In the paper we consider DDGs with these parameters. We prove that if $n$ is odd then such graph can only be a $4\times n$-lattice graph. If $n$ is even we characterise all DDGs with such parameters. Moreover, we characterise all DDGs with parameters $(4n,3n-2,3n-6,2n-2,4,n)$ which are related to $4\times n$-lattice graphs.
\end{abstract}

\begin{keyword} divisible desing graph \sep divisible design \sep Deza graph \sep lattice graph
\vspace{\baselineskip}
\MSC[2010] 05C50\sep 05E10\sep 15A18
\end{keyword}
\end{frontmatter}

\section{Introduction}

A $k$-regular graph is called a {\it divisible design graph} (DDG for short) if its vertex set can be partitioned into $m$ classes of size $n$, such that two distinct vertices from the same class have exactly $\lambda_1$ common neighbors, and two vertices from different classes have exactly $\lambda_2$ common neighbors. A DDG with $m = 1$, $n = 1$, or $\lambda_1 = \lambda_2$ is called improper, in the opposite case it is called proper. For the first time divisible design graphs were studied in master's thesis by M.A.~Meulenberg \cite{M2008} and then they were studied in more detail n two following papers by D.~Crnkovic, W.H.~Haemers, H.~Kharaghani and M.A.~Meulenberg \cite{CH2011, HKM2011} in 2011.

Any graph $G$ can be interpreted as a design, by taking the vertices of $G$ as points, and the neighborhoods of the vertices as blocks. In other words, the adjacency matrix of $G$ is interpreted as the incidence matrix of a design. Let's call this design the {\it neighborhood design} of $G$. An incidence structure with constant block size $k$ is a {\it divisible design} if the point set can be partitioned into $m$ classes of size $n$, such that two points from one class occur together in $\lambda_1$ blocks, and two points from different classes occur together in exactly $\lambda_2$ blocks . A divisible design $D$ is called symmetric (or to have the dual property) if the dual of $D$ (that is the design with the transposed incidence matrix) is again a divisible design with the same parameters as $D$.

The definition of a DDG implies that the neighbourhood design of a DDG is a symmetric divisible design. Conversely, a symmetric divisible design with a polarity with no absolute points is the neighbourhood design of a DDG.  

A {\it Deza graph} with parameters $(v,k,b,a)$ is a $k$-regular graph with $v$ vertices such that any two distinct vertices have $b$ or $a$ common neighbours, where $b \geqslant a$. The definition of a DDG implies that DDG is a Deza graph with $\{b,a\} = \{\lambda_1,\lambda_2\}$.

\medskip

An {\it $m\times n$-lattice graph} is a line  graph of complete bipartite graph $K_{m,n}$. An $m\times n$-lattice graph is a DDG when $m=4$ (or $n=4$). In this paper we characterise DDGs with parameters that are the same as parameters of a of $4\times n$-lattice graph. These graphs have parameters $(4n,n+2,n-2,2,4,n)$. Also we characterise all DDGs with parameters $(4n,3n-2,3n-6,2n-2,4,n)$ which are related to $4\times n$-lattice graphs.

In \cite[Construction 4.8]{HKM2011} the construction of two series of DDGs with parameters $(4n,n+2,n-2,2,4,n)$ and $(4n,3n-2,3n-6,2n-2,4,n)$ from Hadamard matrices was presented. The first one corresponds to $4\times n$-lattice graphs and the second one can be obtained from $4\times n$-lattice graphs by switching all edges between two pairs of rows.
\medskip

The main result of this article is the characterisation of DDGs with parameters $(4n,n+2,n-2,2,4,n)$ and $(4n,3n-2,3n-6,2n-2,4,n)$.
\medskip

The paper is organised as follows. In Section 2 we give some definitions, notations and preliminaries about DDGs and Deza graphs. In section 3 we consider Deza graphs with parameters $(4n,n+2,n-2,2)$ and $(4n,3n-2,3n-6,2n-2)$ which are not DDGs, in particular we prove that in this case $n\leq 8$. In Section 4 we characterise DDGs with parameters $(4n,n+2,n-2,2,4,n)$. In Section 5 we characterise DDGs with parameters $(4n,3n-2,3n-6,2n-2,4,n)$.

\section{Preliminaries}

\subsection{Properties of DDGs}

\begin{proposition}\label{eigenvalues} {\rm  \cite[Lemma 2.1]{HKM2011}}
The eigenvalues of the adjacency matrix of a DDG with parameters $(v,k,\lambda_1,\lambda_2,m,n)$ are
$$\{k, \pm\sqrt{k-\lambda_1},\pm\sqrt{k^2-\lambda_2 v} \},$$ 
with the multiplicities  $1, f_1, f_2, g_1, g_2$, respectively. 
Moreover, $f_1+f_2=m(n-1)$ and $g_1+g_2=m-1$.
\end{proposition}

\medskip

The equation \cite[Equation 2]{HKM2011} states that the trace of adjacency matrix $A$ of a DDG is:

\begin{equation}\label{traceOfA}
    trace(A) = 0 = k + (f_1-f_2)\sqrt{{k-\lambda_1}} + (g_1-g_2)\sqrt{k^2-\lambda_2 v}
\end{equation}

\begin{proposition}\label{equitable}{\rm  \cite[Theorem 3.1]{HKM2011}}
The canonical partition of the adjacency matrix of a proper DDG is equitable, and the quotient matrix $R$ satisfies the following equation

$R^2 = RR^T = (k^2-\lambda_2v)I_m + \lambda_2nJ_m$.
\end{proposition}

\begin{proposition}\label{propertyOfR}{\rm  \cite[Proposition 3.2]{HKM2011}}
The quotient matrix $R$ of DDG satisfies the following conditions

$\sum\limits_i(R)_{i,j} = k$, for all $j = 1,2,\ldots,m$,

$0\leq trace(R) = k+(g_1-g_2)\sqrt{k^2-\lambda_2 v}\leq m(n-1)$.
\end{proposition}

\subsection{Construction which arise from Hadamard matrices}

An $m \times m$ matrix $H$ is a {\it Hadamard matrix} if every entry is $1$ or $-1$, and $H H^T = mI$. A Hadamard matrix $H$ is called {\it graphical} if $H$ is symmetric with constant diagonal, and regular if all row and column sums are equal (for example $l$). Without loss of generality we assume that a graphical Hadamard matrix has diagonal entries $-1$. Let us consider a regular graphical Hadamard matrix $H$. 

\medskip

The next construction is based on \cite[Contruction 4.8]{HKM2011}.

\begin{construction}\label{Hadamard}
Consider the smallest regular graphical Hadamard matrices  \\

    $\left[\begin{array}{cccc}
       -1 & 1 & 1 & 1 \\ 
       1 & -1 & 1 & 1 \\
       1 & 1 & -1 & 1  \\
       1 & 1 & 1 & -1 
\end{array}\right]$ and 
$\left[\begin{array}{cccc}
       -1 & 1 & -1 & -1 \\ 
       1 & -1 & -1 & -1 \\
       -1 & -1 & -1 & 1  \\
       -1 & -1 & 1 & -1 
\end{array}\right].$\\

Replace each entry with value $-1$ by $J_n - I_n$ , and each $+1$ by $I_n$ , then we obtain the adjacency matrix of a DDG with parameters $(4n,n+2,n-2,2,4,n)$ and $(4n,3n-2,3n-6,2n-2,4,n)$. 
\end{construction}

The second graph can be obtained from the first one by switching edges between two pairs of classes of the canonical partition.

\subsection{Switching constructions}

{\it Switching} a set of vertices in a graph means reversing the adjacencies of each pair of vertices, one in the set and the other not in the set: thus, the edge set is changed so that an adjacent pair becomes nonadjacent and a nonadjacent pair becomes adjacent. The edges whose endpoints are both in the set, or both not in the set, are not changed. Graphs are {\it switching equivalent} if one can be obtained from the other by switching. Switching was introduced by van Lint and Seidel (see \cite{LS66}) and developed by Seidel.

An involutive automorphism of a graph is called \emph{Seidel automorphism} if it interchanges only non-adjacent vertices. 

\begin{construction}[{Dual Seidel switching; \cite[Theorem 3.1]{EFHHH1999}}]\label{DSS}  
Let $G$ be a strongly regular graph with parameters $(v,k,\lambda,\mu)$, where $k\neq \mu$, $\lambda \neq \mu$. Let $M$ be the adjacency matrix of $G$, and $P$ be a non-identity permutation matrix of the same size. Then $PM$ is the adjacency matrix of a Deza graph $\Gamma$ if and only if $P$ represents a Seidel automorphism. Moreover, $\Gamma$ is a strictly Deza graph if and only if $\lambda \neq 0$, $\mu \neq 0$.
\end{construction}

\begin{construction}[{Generalised dual Seidel switching 2; \cite[Theorem 6]{KKS2021}}]\label{GDSS2}
Let $G$ be a Deza graph with the adjacency matrix $M$, and $H$ be its induced subgraph with the adjacency matrix $M_{11}$. If there exists a Seidel automorphism of $H$ with the permutation matrix $P_{11}$ such that $P_{11}M_{12}M_{22}=M_{12}M_{22}$, then matrix
$$N= \left(\begin{array}{cc} P_{11}M_{11} & M_{12}\\ M_{21}& M_{22}\\ \end{array}\right)$$
is the adjacency matrix of a Deza graph.
\end{construction}

\begin{remark}\label{rem1} The combinatorial meaning of the matrix condition is as follows. \\
Condition $P_{11}M_{11}M_{12}=M_{11}M_{12}$ from Theorem~\ref{GDSS2} means that for any $v \in V(G)\setminus V(H)$ and for any $x,y \in  V(H)$ such that $\varphi(x)=y$, the number of common neighbours for $v$ and $x$ in $H$ is equal to the number of common neighbours for $v$ and $y$ in $H$.
\end{remark}

\begin{remark}\label{rem4}
Note that in \cite{KKS2021} this Construction was considered only for Deza graphs with strongly regular children but the proof does not use this property, therefore, this Construction can be applied to any Deza graph.
\end{remark}

\section{Deza graphs which are not DDGs}

\begin{proposition}\label{a<2b-k} {\rm  \cite[Theorem 1]{KS2019}}
If $G$ is a Deza graph with parameters $(v,k,b,a)$ and $a<2b-k$, 
then $G$ is a DDG.
\end{proposition}

\begin{lemma}
If $G$ is a Deza graph with parameters $(4n,n+2,n-2,2)$ or $(4n,3n-2,3n-6,2n-2)$ and $G$ is not a DDG, then $n \leq 8$. 
\end{lemma}

\proof 
\begin{enumerate}
    \item Let $G$ have parameters $(4n,n+2,n-2,2)$. By lemma \ref{a<2b-k}, if $G$ is not a DDG, then $a\geq 2b-k$. Thus, $2\geq 2(n-2) - (n+2)$, and $n\leq 8$.
    \item Let $G$ have parameters $(4n,3n-2,3n-6,2n-2)$. By lemma \ref{a<2b-k}, if $G$ is not a DDG, then $a\geq 2b-k$. Thus, $2n-2\geq 2(3n-6) - (3n-2)$, and $n\leq 8$. \hfill $\square$
\end{enumerate}

Deza graphs with parameters $(4n,n+2,n-2,2)$ and $(4n,3n-2,3n-6,2n-2)$ in case $n \leq 8$ were determined completely by computer search. We found $48$ non-isomorphic Deza graphs with parameters $(4n,n+2,n-2,2)$ and $10$ non-isomorphic Deza graphs with parameters $(4n,3n-2,3n-6,2n-2)$ for $n=6$  which are not DDGs. For the remaining $n$, we found only DDGs. 

Adjacency matrices of all Deza graphs with such parameters are available on the web pages \url{http://alg.imm.uran.ru/dezagraphs/deza.php?v=24&k=8&b=4&a=2&form=None} and \url{http://alg.imm.uran.ru/dezagraphs/deza.php?v=24&k=16&b=12&a=10&form=None}.   

Since we found all DDGs up to $32$ vertices, further we assume that $n > 8$.

\section{DDGs with parameters $(4n,n+2,n-2,2,4,n)$}

The main result of this section is the following theorem.

\begin{theorem}\label{mainThm}
Let $G$ be a DDG with parameters $(4n,n+2,n-2,2,4,n)$. 
Then one of the following cases hold.
\begin{enumerate}
    \item If $n$ is odd, then $G$ is isomorphic to $4\times n$-lattice graph,
    \item If $n$ is even, then the quotient matrix $R$ of $G$ equals one of matrices~(\ref{matrix1}), (\ref{matrix2}) or~(\ref{matrix3})
    \begin{enumerate}
        \item if $R$ equals matrix~(\ref{matrix1}), then $G$ is isomorphic to $4\times n$-lattice graph,
        \item if $R$ equals matrix~(\ref{matrix2}), then $G$ is isomorphic to graph $G'$ from Construction~\ref{matrix2GDSS},
        \item if $R$ equals matrix (\ref{matrix3}), then $G$ is isomorphic to one of the graphs obtained from Construction~\ref{reverseSwitching}.
    \end{enumerate}
    
\end{enumerate}
\end{theorem}

By Proposition~\ref{eigenvalues} we have $g_1+g_2 = m-1 = 3$. 
We can calculate all possibilities for $g_1,g_2$ and $tr(R)$ using Proposition~\ref{propertyOfR} and Equation~\ref{traceOfA}.

\begin{table}[ht]
    \centering
    $$\begin{array}{c|c|c}
       g_1 & g_2 & tr(R) \\
       \hline
       3 & 0 & 4n-4 \\
       2 & 1 & 2n  \\
       1 & 2 & 4 \\
       0 & 3 & 8 - 2n 
    \end{array}$$
    \caption{}
    \label{multiplicities}
\end{table}

\medskip

If $g_1 = 0$ and $n > 8$, then $tr(R) < 0$ and it is impossible.

\medskip

If $g_1 = 3$, then $G$ has exactly four eigenvalues $\{k,\pm 2,n-2\}$. The classification of graphs with the smallest eigenvalue $-2$ (see \cite[Section 3.12]{BCN}) implies that $G$ is isomorphic to $4\times n$-lattice graph.

\subsection{Quotient matrices}

\begin{lemma}\label{diagBlocksForOddn}
Let $G$ be a DDG with parameters $(v,k,\lambda_1,\lambda_2,m,n)$ and let $n$ be odd. If $R = [r_{ij}]$ is the quotient matrix of $G$, then $r_{ii}$ is even for all $i = 1,\ldots, m$.
\end{lemma}

\proof Since $r_{ii}$ is the valency of the subgraph induced by the vertices of an $i$-th class of canonical partition, then $r_{ii}$ is even for odd $n$.\hfill $\square$

\medskip

\begin{proposition}
Let $G$ be a DDG with parameters $(4n,n+2,n-2,2,4,n)$ and let $(a,b,c,d)$ be the row of quotient matrix $R$. Then $\{a,b,c,d\} = \{n-1,1,1,1\}$.
\end{proposition}

\proof By Proposition~\ref{propertyOfR} we have equation
\begin{equation} \label{sumOfRow}
    a+b+c+d = n+2
\end{equation} 
and  by Proposition~\ref{equitable} we have the following equation 
\begin{equation}\label{sumOfSquares}
    a^2+b^2+c^2+d^2 = n^2-2n+4
\end{equation} 
At first we note that if there is $x\in \{a,b,c,d\}$, such that $x \geq n$, then equation (\ref{sumOfSquares}) does not hold. Now denote by $x$ the largest element in $\{a,b,c,d\}$ and denote by $y$ the sum of three other elements. Then by (\ref{sumOfRow}) we have $y = n + 2 - x$. By equation (\ref{sumOfSquares}) we have $x^2 + y^2 \geq n^2-2n+4$. Substituting $y$ we get $x^2 + (n+2-x)^2 \geq n^2-2n+4$. Then $x^2 - (n+2)x + 3n \geq 0$ and either $x\leq \cfrac{n+2 - \sqrt{n^2-8n+4}}{2}$ or $x\geq \cfrac{n+2 + \sqrt{n^2-8n+4}}{2}$. 

In the first case, we have $\cfrac{n+2}{4}\leq \cfrac{n+2 - \sqrt{n^2-8n+4}}{2}$ by (\ref{sumOfRow}) and the fact that  $x$ is the largest element in $\{a,b,c,d\}$. Moreover, $2\sqrt{n^2-8n+4} \leq n+2$ then $3n^2-12n+4 \leq 0$, and then we get $n \leq 2$.

In the second case, let us consider $x\leq n-2$. Then $n-2\geq x\geq \cfrac{n+2 - \sqrt{n^2-8n+4}}{2}$. In this case we have $\sqrt{n^2-8n+4} \leq n-6$ and then $n\leq 8$. 

If $x = n-1$ then $y = 3$ and by (\ref{sumOfSquares}) other three elements in $\{a,b,c,d\}$ equal $1$. Since $n>8$ then Proposition 5 is proved. \hfill $\square$

\begin{corollary}\label{quotientMatrices}
Let $G$ be a DDG with parameters $(4n,n+2,n-2,2,4,n)$. Then there are exactly three following possibilities for quotient matrix $R$ of $G$.
\end{corollary}

\begin{equation}\label{matrix1}
    \left[\begin{array}{cccc}
       n-1 & 1 & 1 & 1 \\ 
       1 & n-1 & 1 & 1 \\
       1 & 1 & n-1 & 1 \\
       1 & 1 & 1 & n-1 
\end{array}\right]
\end{equation}

\begin{equation}\label{matrix2}
    \left[\begin{array}{cccc}
       1 & n-1 & 1 & 1 \\ 
       n-1 & 1 & 1 & 1 \\
       1 & 1 & n-1 & 1 \\
       1 & 1 & 1 & n-1 
\end{array}\right]
\end{equation}

\begin{equation}\label{matrix3}
    \left[\begin{array}{cccc}
       1 & n-1 & 1 & 1 \\ 
       n-1 & 1 & 1 & 1 \\
       1 & 1 & 1 & n-1 \\
       1 & 1 & n-1 & 1 
\end{array}\right]
\end{equation}

\medskip

\begin{remark}\label{remMatrix1}
It is easy to see that these matrices correspond to rows of table \ref{multiplicities}. Therefore, if graph $G$ has quotient matrix (\ref{matrix1}) then $G$ is isomorphic to $4\times n$-lattice graph.
\end{remark}

{\bf Statement 1 of Theorem \ref{mainThm}} immediately follows from Lemma~\ref{diagBlocksForOddn} and Remark~\ref{remMatrix1}.
\medskip

Further we will assume that $n$ is even. We need to consider graphs with quotient matrices (\ref{matrix2}) and (\ref{matrix3}).

\subsection{Graphs with quotient matrix~(\ref{matrix2})}

\begin{construction}\label{matrix2GDSS}
The $4\times n$-lattice graph $G$ has an induced subgraph $H$ such that $H$ is isomorphic to $2\times n$-lattice graph. By Remark~\ref{rem4}, Theorem~\ref{GDSS2} holds for $G$ and $H$. For any even $n \geqslant 6$, there is a Seidel automorphism $\varphi$ of $H$ corresponding to the central symmetry of the lattice with two rows and $n$ columns. The subgraph $H$ satisfies the condition of Theorem~\ref{GDSS2} since the combinatorial condition of Remark~\ref{rem1} holds. Indeed, for any $v \in V(G)\setminus V(H)$ and for any $x,y \in V(H)$ such that $\varphi(x)=y$, there is only one neighbour in $H$ for both pairs $v$ and $x$, and $v$ and $y$. Hence, by Theorem~\ref{GDSS2}, for $n \geqslant 6$ there is a Deza graph with parameters $(n^2,2(n-1),n-2,2)$. Since automorphism $\varphi$ interchanges two blocks of canonical partition of $4\times n$-lattice graph, obtained Deza graph is a DDG. And it is clear that it has quotient matrix (\ref{matrix2}). Denote this graph by $G'(n)$.
\end{construction}

In this section we prove the statement $2b$ of Theorem \ref{mainThm}. If $G$ is a DDG with parameters $(4n,n+2,n-2,2,4,n)$ and  quotient matrix (\ref{matrix2}), then $G$ is isomorphic to $G'(n)$.

\proof Denote the blocks of the canonical partition of $G$ by $V_1$, $V_2$, $V_3$ and $V_4$ with respect to the quotient matrix~(\ref{matrix2}). According to the quotient matrix~(\ref{matrix2}), the vertices of $V_1$ and $V_2$ induce cliques in $G$. There is perfect matching between them. We denote the vertices of the first block by $(1,1), (1,2), \ldots, (1, n)$, and the vertices of the second block by $(2,1), (2,2), \ldots, ( 2, n)$, such that vertices $(1, i)$ and $(2, i)$ are adjacent for any $i$. 

Consider all pairs of vertices from blocks $V_1$ and $V_2$. 
Recall that in $G$ vertices from the same block have $n-2$ common neighbours, and vertices from the different blocks have $2$ common neighbours. Any two vertices from $V_1$ have $n-2$ common neighbours in $V_1$, and no common neighbours in other blocks. It is also true for any two vertices from $V_2$. Vertices $(1,i)$ and $(2,j)$, where $i\neq j$, have two common neighbours ($(1,j)$ and $(2,i)$) in $V_1\cup V_2$, and no common neighbours in $V_3$ and $V_4$. Vertices $(1,i)$ and $(2,i)$ have no common neighbours in $V_1\cup V_2$, and $2$ common neighbours in $V_3\cup V_4$. But $(1,i)$ has only one neighbour in $V_3$ and one neighbour in $V_4$. Denote these two vertices by $(3,i)$ in $V_3$ and $(4,i)$ in $V_4$. Then $(3,i)$ and $(4,i)$ are adjacent with the vertex $(2,i)$ and have no other neighbours in $V_1 \cup V_2$. 

Now consider common neighbours of vertices $(1,i)$ and $(3,i)$. In this case $N((1,i))\cap N((3,i)) = \{(2,i),(4,i)\}$ because $N((1,i))\cap (V_3\cup V_4) = \{(3,i),(4,i)\}$, $N((3,i))\cap (V_1\cup V_2) = \{(1,i),(2,i)\}$ and vertices $(1,i)$, $(3,i)$ have $2$ common neighbours in graph $G$. Similarly, $N((1,i))\cap N((4,i)) = \{(2,i),(3,i)\}$ hence, vertices $\{(1,i),(2,i),(3,i),(4,i)\}$ induce a clique in graph $G$.

From the quotient matrix~(\ref{matrix2}) we know that the subgraph of $G$ induced by $V_3$ is isomorphic to $n/2$ copies of $K_2$. Without loss of generality we assume that any vertex $(3,i)$ is adjacent to the vertex $(3,n-i+1)$ for all $i=1,2,\ldots, n/2$. Now consider vertices $(4,i)$ and $(4,n-i+1)$. Since vertices $(1,i)$ and $(3,n-i+1)$ have $2$ common neighbours in $G$ and $(1,n-i+1),(3,i) \in N((1,i),(3,n-i+1))$, then any vertex $(3,n-i+1)$ is not adjacent to the vertex $(4,i)$. Similarly, any vertex $(3,i)$ is not adjacent to the vertex $(4,n-i+1)$. Therefore, vertices $(1,i)$ and $(4,n-i+1)$ have only one common neighbour in $V_1\cup V_2 \cup V_3$. But any vertex $(1,i)$ has one neighbour (the vertex $(4,i)$) in $V_4$. Hence, vertices $(4,i)$ and $(4,n-i+1)$ are adjacent for any $i$. 

We described all edges of the graph $G$. Hence, there is the unique DDG with such parameters and quotient matrix~(\ref{matrix2}). Thus, this graph is isomorphic to $G'(n)$.\hfill $\square$

\subsection{Graphs with quotient matrix~(\ref{matrix3})}

Let $G$ be a DDG with parameters $(4n,n+2,n-2,2,4,n)$ and the quotient matrix~(\ref{matrix3}). 
Further, let $A = [A_{ij}]$ be the adjacency matrix of $G$ with blocks $A_{ij}$ corresponding to a canonical partition of $G$. Now consider an auxiliary graph $G^\ast$ with the following adjacency matrix

\begin{equation}\label{matrix4}
   A^\ast =  \left[\begin{array}{cccc}
       A_{11} & J_n - A_{12} & A_{13} & A_{14} \\ 
       J_n - A_{21} & A_{22} & A_{23} & A_{24} \\
       A_{31} & A_{32} & A_{33} & J_n - A_{34} \\
       A_{41} & A_{42} & J_n - A_{43} & A_{44} 
\end{array}\right].
\end{equation}

This matrix can be obtained by the operation that corresponds to the switching of all edges between $V_1$ and $V_2$, and between $V_3$ and $V_4$. 
Graph $G^\ast$ with the adjacency matrix $A^\ast$ has an equitable partition with the quotient matrix $J_4$. We assume that the sets of vertices $G$ and $G^\ast$ coincide and the canonical partitions of graph $G$ are the same as the equitable partition of $G^\ast$.

Denote classes of the equitable partitions in graphs $G$ and $G^\ast$ by $V_1$, $V_2$, $V_3$ and $V_4$. Note that graph $G^\ast$ is regular with valency $4$ and it can be disconnected. In the following lemma we calculate the number of common neighbours of all pairs of vertices in graph $G^\ast$. Denote the set of common neighbours of vertices $x$ and $y$ in graph $G^\ast$ by $N_{G^\ast}(x,y)$. Each vertex $x\in V_i$ in graphs $G$ and $G^\ast$ has the only one neighbour in its own class $V_i$ and we denote this neighbour by~$x'$.

\begin{lemma} \label{pairsOfVertices}
Let $G^\ast$, $V_1$, $V_2$, $V_3$ and $V_4$ be as described above, and $x$ and $y$ be distinct vertices of $G^\ast$. 
Without loss of generality we can assume that $x\in V_1$. Then one of the following cases hold:
\begin{enumerate}
    \item If $y\in V_1$, then $|N_{G^\ast}(x,y)| = 0$,
    \item Let $y\in V_2$, and
    \begin{enumerate}
        \item if $|N(x,y) \cap (V_1\cup V_2)| = 2$, then $|N_{G^\ast}(x,y)| = 0$,
        \item if $|N(x,y) \cap (V_1\cup V_2)| = 1$, then $|N_{G^\ast}(x,y)| = 2$,
        \item if $|N(x,y) \cap (V_1\cup V_2)| = 0$, then $|N_{G^\ast}(x,y)| = 4$,
    \end{enumerate}
    \item Let $y\in V_i$, where $i\in \{3,4\}$, and
    \begin{enumerate}
        \item if $|N(x,y) \cap (V_1\cup V_i)| = 2$, then $|N_{G^\ast}(x,y)| = 4$,
        \item if $|N(x,y) \cap (V_1\cup V_i)| = 1$, then $|N_{G^\ast}(x,y)| = 2$,
        \item if $|N(x,y) \cap (V_1\cup V_i)| = 0$, then $|N_{G^\ast}(x,y)| = 0$,
    \end{enumerate}
\end{enumerate}
\end{lemma}

\proof
\begin{enumerate}
    \item Since $G$ has qoutient matrix \ref{matrix3}, then all common neighbours of $x$ and $y$ in $G$ lie in $V_2$. Since the edges between $V_1$ and $V_2$ in $G^\ast$ form a perfect matching then $x$ and $y$ have no common neighbours in $G^\ast$,
    \item 
    \begin{enumerate}
        \item The vertex $x$ is adjacent to $y'$ and the vertex $y$ is adjacent to $x'$ in $G$, but in $G^\ast$ these adjacencies are removed. Thus, $x$ and $y$ have no common neighbours in $G^\ast$.
        \item In $G$, either $x$ is adjacent to $y'$ or $y$ is adjacent to $x'$ and they have one more common neighbour in $V_3\cup V_4$. In $G^\ast$, edges between $V_1$ and $V_2$ are switched then $x$ and $y$ have one common neighbour in $V_1\cup V_2$ and one common neighbour in $V_3\cup V_4$.
        \item In $G$, neither $x$ is adjacent to $y'$ nor $y$ is adjacent to $x'$. Hence $x$ is adjacent to $y'$ and $y$ is adjacent to $x'$ in $G^\ast$. Moreover $x$ and $y$ have two common neighbours in $V_3\cup V_4$.
    \end{enumerate}
    \item Let $y\in V_3$. If $y\in V_4$, then the proof is the same as in the previous point. 
    \begin{enumerate}
        \item In $G$, vertices $x$ and $y$ have two common neighbours in $V_1\cup V_3$ and 
        have no common neighbours in $V_2\cup V_4$. Then, in $G^\ast$, they have two common neighbours in $V_1\cup V_3$ and two common neighbours in $V_2\cup V_4$.
        \item In $G$, vertices $x$ and $y$ have one common neighbour in $V_1\cup V_3$ and one common neighbour in $V_2\cup V_4$. Then, in $G^\ast$, they have the same one common neighbour in $V_1\cup V_3$ and the other common neighbour in $V_2\cup V_4$.
        \item In $G$, vertices $x$ and $y$ have no common neighbours in $V_1\cup V_3$ and two common neighbours in $V_2\cup V_4$. Then, in $G^\ast$, they have no common neighbours in $V_1\cup V_3$ and in $V_2\cup V_4$.
    \end{enumerate}
\end{enumerate} \hfill $\square$

Consider a connected component of $G^\ast$. Denote it by $D$ and let $D_1 = D\cap V_1$, $D_2 = D\cap V_2$, $D_3 = D\cap V_3$, $D_3 = D\cap V_3$. Since $G^\ast$ has the equitable partition with parts $V_1$, $V_2$, $V_3$, $V_4$ and with quotient matrix $J_4$, then the set $D_i\cup D_j$ induce subgraph of valency $2$, which is the union of cycles, for any $i$ and $j$ ($i\neq j$).

\begin{lemma}\label{sizeDevidedBy8}
The size of $D$ is divisible by $8$.
\end{lemma}

\proof Consider the largest independent set $S$ in $D_1$. For each vertex $x\in S$ we have one neighbour $x'$ in $D_1$. Also for each vertex $x$ we have one neighbour in $D_2$, one neighbour in $D_3$ and one neighbour in $D_4$. It is also true for the vertex $x'$. Since any other vertex must have one neighbour in $D_1$ which is not adjacent with vertices from $S$, then there are no other vertices in $D$. Finally, we have $|D| = 8|S|$ and $|D|$ is divisible by $8$.\hfill $\square$ 

\medskip

Now we assume that $D$ have vertices $x$ and $y$ with $4$ common neighbors.

\begin{lemma}\label{2c}
There exist vertices $x,y \in D$ of type $2c$ from Lemma~\ref{pairsOfVertices} if and only if there are four vertices in $D_1\cup D_2$ (or in $D_3\cup D_4$) which induce a cycle. 
In this case, the induced subgraphs on $D_1\cup D_2$ and $D_3\cup D_4$ are isomorphic to $s C_4$ for some $s$. Moreover, $D$ is isomorphic to $C_s[\overline{K_2}]$ where each copy of $\overline{K_2}$ is a pair of vertices of type $2c$ from Lemma~\ref{pairsOfVertices}. Alternately, two such pairs are from $V_1\cup V_2$ and two pairs are from $V_3\cup V_4$. 
\end{lemma}

\proof At first consider vertices $x,y\in D$ of type $2c$ from Lemma~\ref{pairsOfVertices}. They have two common neighbours ($x'$ and $y'$) in $D_1\cup D_2$ and two common neighbours (say $z$ and $t$) in $D_3\cup D_4$. Hence the pairs $x', y'$ and $z, t$ are of type $2c$. Then $ z', t' $ and two common neighbors $ x', y' $ in $ D_3 \cup D_4 $ are of type $ 2c$. Then the common neighbors of $ z', t' $ are also of type $ 2c $. Let's continue these process until the cycle of pairs of vertices is closed. Thus, we have a cycle of pairs of type $2c$, where two pairs $ x, y $ and $ x', y' $ are taken from $ D_1 \cup D_2 $, the next two pairs are from $ D_3 \cup D_4 $ and further alternately, because we can start with any of these pairs.

Conversely, if we have the cycle $(x,x',y,y')$ as the induced subgraph in $D_1\cup D_2$ (or $D_3\cup D_4$) then the pairs of vertices $x,y$ and $x',y'$ are of type $2c$, and the lemma is proved.\hfill $\square$

\begin{lemma}\label{3a}
There exist vertices $x,y \in D$ of type $3a$ from Lemma~\ref{pairsOfVertices} if and only if there are four vertices in $D_1\cup D_3$ (or in $D_2\cup D_4$) which induce cycle. 
In this case, the induced subgraphs on $D_1\cup D_2$ and $D_3\cup D_4$ are isomorphic to $s C_4$ for some $s$. Moreover, $D$ is isomorphic to $C_s[\overline{K_2}]$ where each copy of $\overline{K_2}$ is a pair of vertices of type $2c$ from Lemma~\ref{pairsOfVertices}. Alternately, two such pairs are from the set $V_1\cup V_3$ and two pairs are from $V_2\cup V_4$. 
\end{lemma}

\proof At first consider vertices $x,y\in D$ of type $3a$ from Lemma~\ref{pairsOfVertices}. They have two common neighbours ($x'$ and $y'$) in $D_1\cup D_3$ and two common neighbours (say $z$ and $t$) in $D_2\cup D_4$. Then pairs $x', y'$ and $z, t$ are of type $3a$. Then $ z', t' $ and two common neighbors $ x', y' $ in $ D_2 \cup D_4 $ are of type $ 3a $. Then the common neighbors of $ z', t' $ are also of type $ 3a $. Let's continue these process until the cycle of pairs of vertices is closed. Thus, we have a cycle of pairs of type $ 3a $, where two pairs $ x, y $ and $ x', y' $ are taken from $ D_1 \cup D_3 $, the next two pairs are from $ D_2 \cup D_4 $ and further alternately, because we can start with any of these pairs.

Conversely, if we have cycle $(x,x',y,y')$ as induced subgraph in $D_1\cup D_3$ (or $D_2\cup D_4$) then the pairs of vertices $x,y$ and $x',y'$ are of type $3a$, and the lemma is proved.\hfill $\square$

\medskip

Subgraph $D$ in Lemma~\ref{2c} and Lemma~\ref{3a} is the same one but it is embedded in a different way in graph $G^\ast$. Hence, corresponding subgraphs in graph $G$ are not isomorphic (excluding the case when for $D$ both lemmas holds). 

\begin{lemma}\label{noPairWith4Neibours}
If connected component $D$ has no pairs of vertices with four common neighbours, then $D$ is isomorphic to $4$-cube.
\end{lemma}

\proof If connected component $D$ has no pairs of vertices with four common neighbours then $D$ is a $(0,2)$-graph of valency $4$. Mulder in \cite{M79} proved that if graph is a $k$-regular $(0,2)$-graph on $v$ vertices and has the diameter $d$, then $v\leq 2^k$ and $d\leq k$. In both cases equality  is true only  for the $k$-cube  with parameters $(2^k, k, 2, 0)$. By Lemma~\ref{sizeDevidedBy8} the size of $D$ is divisible by $8$, so $D$ is isomorphic to $4$-cube or is a Deza graph with parameters $(8,4,2,0)$. But in the last case  $D_1\cup D_2$ induces cycle $C_4$, so by Lemma~\ref{2c} we have a contradiction. \hfill $\square$

\begin{lemma}
There exist three non-isomorhic equitable partitions of $4$-cube with quotient matrix $J_4$.
\end{lemma}

\proof The inner edges of the parts of the partition form a perfect matching in $4$-cube. Moreover, two edges from the same part are antipodal in $4$-cube. Now we can use the result that was obtained in paper \cite{GH88} that states that there are eight equivalence classes of perfect matchings in $4$-cube and only three of them satisfy necessary condition. It is enough to test these perfect matchings to prove that there are three antipodal perfect matchings and three corresponding equitable partitions.\hfill $\square$

\begin{construction}\label{reverseSwitching}
Consider some copies of graph $C_s[\overline{K_2}]$ with equitable partitions as in Lemmas \ref{2c} or \ref{3a} and some copies of $4$-cube with one of the three equitable partitions. Then we can switch edges between $V_1$ and $V_2$ and also between $V_3$ and $V_4$. This switching gives us the DDG with parameters $(4n,n+2,n-2,2,4,n)$ and quotient matrix (\ref{matrix3}). 
\end{construction}

To complete the proof of Theorem \ref{mainThm}, we prove statement $2c$.\medskip

{\bf Proof of statement 2c of Theorem \ref{mainThm}}. By Lemmas~\ref{2c}, \ref{3a} and \ref{noPairWith4Neibours} the connected component $D$ of graph $G^\ast$ is isomorphic to $C_s[\overline{K_2}]$ or $4$-cube. Graph $C_s[\overline{K_2}]$ has two non-isomorphic embeddings into  canonical partition of graph $G^\ast$. $4$-cube has three non-isomorphic embeddings into  canonical partition of graph $G^\ast$. Now we can obtain graph $G^\ast$ only from such connected components. Hence divisible design graph $G$ can be obtained from $G^\ast$ by the reverse switching. Thus, Theorem~\ref{mainThm} is proved. \hfill $\square$

\section{DDGs with parameters $(4n,3n-2,3n-6,2n-2,4,n)$}

The main result of this section is the following theorem.

\begin{theorem}\label{mainThm2}
Let $G$ be a DDG with parameters $(4n,3n-2,3n-6,2n-2,4,n)$. Then $G$ can be obtained from a DDG with parameters $(4n,n+2,n-2,2,4,n)$ by switching between the first two classes and the last two classes of the canonical partition of DDG. Classes are numbered according to the quotient matrix from corollary \ref{quotientMatrices}.
\end{theorem}

\subsection{Quotient matrices}
 
By Proposition~\ref{eigenvalues} we have $g_1+g_2 = m-1 = 3$. We can calculate all possibilities for $g_1,g_2,f_1,f_2$ and $tr(R)$ using Proposition~\ref{propertyOfR} and Equation~\ref{traceOfA}.

\begin{table}[ht]
    \centering
    $$\begin{array}{c|c|c}
       g_1 & g_2 & tr(R)\\
       \hline
       3 & 0 & 6n-8\\
       2 & 1 & 4n-4 \\
       1 & 2 & 2n \\
       0 & 3 & 4
    \end{array}$$
    \caption{}
    \label{multiplicities2}
\end{table}

\medskip

If $g_1 = 3$ and $n > 8$, then $tr(R) > 4n-4$ and it is impossible by proposition \ref{propertyOfR} 

\medskip

\begin{proposition}
Let $G$ be a DDG with parameters $(4n,3n-2,3n-6,2n-2,4,n)$ and let $(a,b,c,d)$ be the row of the quotient matrix $R$. Then $\{a,b,c,d\} = \{n-1,n-1,n-1,1\}$.
\end{proposition}

\proof By Proposition~\ref{propertyOfR} we have the equation
\begin{equation} 
    a+b+c+d = 3n-2
\end{equation} 
and  by Proposition~\ref{equitable} we have the following equation 

\begin{equation}
    a^2+b^2+c^2+d^2 = 3n^2-6n+4
\end{equation} 
Assume that $d$ is the smallest element in $\{a,b,c,d\}$. At first we note that $a^2+b^2+c^2+d^2$ is maximal if $a$, $b$ and $c$ are maximal. If $d \geq 2$ then $a^2 + b^2 + c^2$ is maximal when $a=n$, $b = n$, $c = n-4$, and $d = 2$. But in this case $a^2 + b^2 + c^2 +d^2 = 3n^2 - 8n + 20$. While $n>8$ this sum of squares is smaller than $3n^2-6n+4$. Now assume that $d = 1$, then $\{a,b,c\}$ is equal to $\{n,n,n-3\}$ or $\{n,n-1,n-2\}$ or $\{n-1,n-1,n-1\}$. Only the last case gives us equality $a^2 + b^2 + c^2 + 1 = 3n^2 + 6n + 4$. If $d = 0$, then $\{a,b,c\}$ is equal to $\{n,n,n-2\}$ or $\{n,n-1,n-1\}$. Both cases give us inequality $a^2 + b^2 + c^2 \neq 3n^2 + 6n + 4$. \hfill $\square$

\begin{corollary}\label{quotientMatrices2}
Let $G$ be a DDG with parameters $(4n,3n-2,3n-6,2n-2,4,n)$. Then there exist exactly 
three following possibilities for quotient matrix $R$ of $G$.
\end{corollary}

\begin{equation}\label{matrix21}
    \left[\begin{array}{cccc}
       n-1 & 1 & n-1 & n-1 \\ 
       1 & n-1 & n-1 & n-1 \\
       n-1 & n-1 & n-1 & 1 \\
       n-1 & n-1 & 1 & n-1 
\end{array}\right]
\end{equation}

\begin{equation}\label{matrix22}
    \left[\begin{array}{cccc}
       n-1 & 1 & n-1 & n-1 \\ 
       1 & n-1 & n-1 & n-1 \\
       n-1 & n-1 & 1 & n-1 \\
       n-1 & n-1 & n-1 & 1 
\end{array}\right]
\end{equation}

\begin{equation}\label{matrix23}
    \left[\begin{array}{cccc}
       1 & n-1 & n-1 & n-1 \\ 
       n-1 & 1 & n-1 & n-1 \\
       n-1 & n-1 & 1 & n-1 \\
       n-1 & n-1 & n-1 & 1 
\end{array}\right]
\end{equation}

\begin{remark}\label{remMatrix21}
These matrices correspond to the rows of table \ref{multiplicities2}. Therefore, if graph $G$ has quotient matrix (\ref{matrix21}) then $G$ is isomorphic to the second graph from construction \ref{Hadamard}.
\end{remark}

If we switch edges between the first two and the second two classes of canonical partition of graphs with quotient matrices (\ref{matrix22}) and (\ref{matrix23}), then we obtain graphs with corresponding equitable partition with quotient matrices (\ref{matrix2}) and (\ref{matrix3}) respectively. By Theorem \ref{mainThm2} it is enough to prove that the resulting graphs are DDGs with parameters $(4n,n+2,n-2,2,4,n)$.

\subsection{Proof of Theorem \ref{mainThm2}}

The proof is carried out by a simple check of the numbers of common neighbours of the pairs of vertices in all possible cases.

Denote the blocks of the canonical partition of $G$ by $V_1$, $V_2$, $V_3$ and $V_4$ with respect to the quotient matrix.

\begin{lemma}\label{lemma9}
Let $G$ be a DDG with parameters $(4n,3n-2,3n-6,2n-2,4,n)$ with quotient matrix (\ref{matrix22}). Let $G'$ be a graph obtained from $G$ by switching between the first two and the last two classes of the canonical partition. Then $G'$ is a DDG with parameters $(4n,n+2,n-2,2,4,n)$ and with quotient matrix (\ref{matrix2}).
\end{lemma}

\proof We need to consider all pairs of vertices in $G$ and show that two vertices from the same class have $n-2$ common neighbours in $G'$ and two vertices from the different blocks have $2$ common neighbours in $G'$. Let's consider all possibilities for number of common neighbours of vertices $x,y\in G$.
\begin{itemize}
    \item Consider $x,y \in V_i$, where $i = 1$ or $2$. In this case $x$ and $y$ have $n-2$ common neighbours in $V_i$, $n-2$ common neighbours in $V_3$ and $n-2$ common neighbours in $V_4$. Then each vertex in $V_3\cup V_4$ is adjacent with $x$ or $y$. Hence vertices $x$ and $y$ have only $n-2$ common neighbours in $V_i$ in graph $G'$. 
    \item Consider $x \in V_1$, $y\in V_2$, where $x\sim y$. In this case $x$ and $y$ have no common neighbours in $V_1\cup V_2$. Hence $x$ and $y$ have $n-1$ common neighbours in $V_3$ and $n-1$ common neighbours in $V_4$. Then there are two vertices in $V_3\cup V_4$ that are nonadjacent with $x$ and $y$. Hence vertices $x$ and $y$ have $2$ common neighbours in $V_3\cup V_4$ in graph $G'$. 
    \item Consider $x \in V_1$, $y\in V_2$, where $x\not \sim y$. In this case $x$ and $y$ have $2$ common neighbours in $V_1\cup V_2$. Hence $x$ and $y$ have $n-2$ common neighbours in $V_3$ and $n-2$ common neighbours in $V_4$. Then each vertex in $V_3\cup V_4$ is adjacent with $x$ or $y$. Hence vertices $x$ and $y$ have $2$ common neighbours in $V_1\cup V_2$ in graph $G'$. 
    \item Consider $x,y \in V_i$, where $i = 3$ or $4$. In this case $x$ and $y$ have $n-2$ common neighbours in $V_{7-i}$, $n-2$ common neighbours in $V_1$ and $n-2$ common neighbours in $V_2$. Then each vertex in $V_3\cup V_4$ is adjacent with $x$ or $y$. Hence vertices $x$ and $y$ have only $n-2$ common neighbours in $V_{7-i}$  in graph $G'$. 
    \item Consider $x \in V_3$, $y\in V_4$, where $|N(x,y)\cap (V_3\cup V_4)| = 2$. In this case $x$ and $y$ have $2n-4$ common neighbours in $V_1\cup V_2$. Hence $x$ and $y$ have $n-2$ common neighbours in $V_1$ and $n-2$ common neighbours in $V_2$. Then each vertex in $V_1\cup V_2$ is adjacent with $x$ or $y$. Hence vertices $x$ and $y$ have $2$ common neighbours in $V_3\cup V_4$ in graph $G'$.
    \item Consider $x \in V_3$, $y\in V_4$, where $|N(x,y)\cap (V_3\cup V_4)| = 1$. In this case $x$ and $y$ have $2n-3$ common neighbours in $V_1\cup V_2$. Hence $x$ and $y$ have $n-1$ common neighbours in $V_1$ and $n-2$ common neighbours in $V_2$ or vice versa. Then there are one vertex in $V_1\cup V_2$ that is nonadjacent with $x$ and $y$. Hence vertices $x$ and $y$ have one common neighbour in $V_1\cup V_2$ and one common neighbour in $V_3\cup V_4$ in graph $G'$.
    \item Consider $x \in V_3$, $y\in V_4$ where $|N(x,y)\cap (V_3\cup V_4)| = 0$. In this case $x$ and $y$ have $2n-2$ common neighbours in $V_1\cup V_2$. Hence $x$ and $y$ have $n-1$ common neighbours in $V_1$ and $n-1$ common neighbours in $V_2$. Then there are two vertices in $V_1\cup V_2$ that are nonadjacent with $x$ and $y$. Hence vertices $x$ and $y$ have two common neighbours in $V_1\cup V_2$ in graph $G'$.
\end{itemize} 
Thus, in all cases two vertices from the same class of canonical partition have $n-2$ common neighbours in $G'$ and two vertices from different classes have $2$ common neighbours. Hence, $G'$ is a DDG with parameters $(4n,n+2,n-2,2,4,n)$. \hfill $\square$

\begin{lemma}
Let $G$ be a DDG with parameters $(4n,3n-2,3n-6,2n-2,4,n)$ and with quotient matrix (\ref{matrix23}). Let $G'$ be a graph obtained from $G$ by switching between the first two and the last two classes of the canonical partition. Then $G'$ is a DDG with parameters $(4n,n+2,n-2,2,4,n)$ and with quotient matrix (\ref{matrix3}).
\end{lemma}

\proof The proof is similar to the proof of lemma \ref{lemma9}. We need to consider all pairs of vertices in $G$ and show that two vertices from the same class have $n-2$ common neighbours in $G'$ and two vertices from the different blocks have $2$ common neighbours in $G'$. Let's consider all possibilities for the number of common neighbours of vertices $x,y\in G$. Without loss of generality we can assume that $x\in V_1$.
\begin{itemize}
    \item Consider $y\in V_1$. In this case $x$ and $y$ have $n-2$ common neighbours in $V_2$, $n-2$ common neighbours in $V_3$ and $n-2$ common neighbours in $V_4$. Then each vertex in $V_3\cup V_4$ adjacent with $x$ or $y$. Hence vertices $x$ and $y$ have only $n-2$ common neighbours in $V_2$ in graph $G'$. 
    \item Consider $y \in V_i$ where $i\in \{2,3,4\}$. Let $\{j,s\} = \{2,3,4\}\setminus \{i\}$.
    \begin{itemize}
        \item If $|N(x,y)\cap (V_1\cup V_i)| = 2$. Then $x$ and $y$ have $2n-4$ common neighbours in $V_j\cup V_s$. Hence $x$ and $y$ have $n-2$ common neighbours in $V_j$ and $n-2$ common neighbours in $V_s$. Then each vertex in $V_j\cup V_s$ is adjacent with $x$ or $y$. Hence vertices $x$ and $y$ have two common neighbours in $V_1\cup V_i$ in graph $G'$.
        \item If $|N(x,y)\cap (V_1\cup V_i)| = 1$. Then $x$ and $y$ have $2n-3$ common neighbours in $V_j\cup V_s$. Hence $x$ and $y$ have $n-1$ common neighbours in $V_j$ and $n-2$ common neighbours in $V_s$ or vice versa. Then there is one vertex in $V_j\cup V_s$ that is nonadjacent with $x$ and $y$. Hence $x$ and $y$ have one common neighbour in $V_1\cup V_i$ and one common neighbour in $V_j\cup V_s$ in graph $G'$.
        \item If $|N(x,y)\cap (V_1\cup V_i)| = 0$. Then $x$ and $y$ have $2n-2$ common neighbours in $V_j\cup V_s$. Hence $x$ and $y$ have $n-1$ common neighbours in $V_j$ and $n-1$ common neighbours in $V_s$. Then there are two vertices in $V_j\cup V_s$ that are nonadjacent with $x$ and $y$. Hence vertices $x$ and $y$ have two common neighbours in $V_j\cup V_s$ in graph $G'$. 
    \end{itemize}
\end{itemize} 
Thus, in all cases two vertices from the same class of canonical partition have $n-2$ common neighbours in $G'$ and two vertices from different classes have $2$ common neighbours. Hence, $G'$ is a DDG with parameters $(4n,n+2,n-2,2,4,n)$. \hfill $\square$

\medskip

This completes the proof of Theorem~\ref{mainThm2}.\hfill $\square$

\begin{remark}
Graphs with quotient matrix \ref{matrix23} have four eigenvalues $\{3n-2, \pm 2, -(n-2)\}$. Hence, these graphs are walk regular.
\end{remark}

\section*{Acknowledgements} 
Author is supported by RFBR according to the research project 20-51-53023.

\end{document}